\documentclass[12pt]{amsart}
\usepackage{amssymb}
\usepackage{pstricks}

\newtheorem{theorem}{Theorem}
\newtheorem{proposition}{Proposition}
\newtheorem{lemma}[proposition]{Lemma}
\newtheorem{corollary}[proposition]{Corollary}

\theoremstyle{definition}
\newtheorem*{definition}{Definition}
\newtheorem*{remark}{Remark}

\newcommand{\Z}{\mathbb{Z}}
\newcommand{\R}{\mathbb{R}}
\newcommand{\C}{\mathbb{C}}
\newcommand{\CP}{\mathbb{CP}}

\newcommand{\dbar}{\bar\partial}

\title{Infinitely many monotone Lagrangian tori in $\R^6$}
\author{Denis Auroux}
\address{Department of Mathematics, UC Berkeley, Berkeley CA 94720-3840, USA}
\email{auroux@math.berkeley.edu}
\thanks{This work was partially supported by NSF grants DMS-1406274 and
DMS-1264662, and by a fellowship from the Simons Foundation.}

\begin{document}
\begin{abstract} 
We construct infinitely many families of monotone Lagrangian tori in
$\R^6$, no two of which are related by Hamiltonian isotopies (or
symplectomorphisms). These families
are distinguished by the (arbitrarily large) numbers of families of 
Maslov index 2 pseudo-holomorphic discs that they bound.
\end{abstract}

\maketitle

\section{Introduction}
The study and classification of Lagrangian submanifolds in symplectic
manifolds is a central topic of modern symplectic topology; in spite of
spectacular advances in the last few decades, it remains poorly
understood, even in very simple symplectic manifolds such as the standard
symplectic vector space $(\R^{2d},\omega_0)$.

By a celebrated result of Gromov, there are no
closed exact Lagrangian submanifolds in $\R^{2d}$, and in fact
any closed Lagrangian in $\R^{2d}$ must bound some 
pseudo-holomorphic discs of non-zero area \cite{Gromov}. 
(This is in sharp contrast with the situation for immersed Lagrangians, 
see e.g.\ \cite{EEMS}.)
Thus, the nicest condition that one could impose on a closed Lagrangian
submanifold $L\subset \R^{2d}$ is for it to be {\em monotone},
i.e.\ that the symplectic area of discs with boundary on $L$ is
(positively) proportional to their {\em Maslov index}.

The simplest examples of monotone Lagrangians in $\R^{2d}$ are the
tori obtained as products of $d$ circles of equal radius, $L=S^1(r)\times\dots \times
S^1(r)$. In the early 1990s Chekanov found the first examples of Lagrangian
tori in $\R^{2d}$ that cannot be related to product tori by Hamiltonian
isotopies (or symplectomorphisms) \cite{Chekanov} (see also \cite{EP}). Subsequent work of Chekanov and Schlenk has led to more
examples, the so-called {\em monotone twist tori} \cite{Chekanov-Schlenk};
the number of tori produced by this construction grows exponentially with 
the dimension, but remains finite for all $d$.

More recently, Renato Vianna's thesis \cite{Vianna} shows that $\CP^2$
contains at least one new kind of monotone Lagrangian torus besides
product and Chekanov tori; this result was recently improved to show that
$\CP^2$ contains {\em infinitely many} non-isotopic monotone
Lagrangian tori \cite{GM,Vianna2}.

In this paper, we construct infinitely many families of monotone Lagrangian
tori in $\R^6$, no two of which are related by symplectomorphisms.
Specifically, the invariants that
we use to distinguish these tori are the algebraic counts of
Maslov index 2 pseudo-holomorphic discs whose boundary passes through a
given point (see \S \ref{ss:count}); these invariants were already used
by Eliashberg-Polterovich to distinguish the Chekanov torus in $\R^4$
\cite{EP} and in much of the subsequent work
\cite{Chekanov-Schlenk,Vianna,Vianna2}.

\begin{theorem}\label{thm:main}
For each integer $n\ge 0$, and for any choice of monotonicity constant,
there exists a monotone Lagrangian torus 
$L\subset (\R^6,\omega_0)$ such that there are $n+2$ distinct 
Maslov index $2$ classes in $\pi_2(\R^6,L)$ for which the algebraic
count of pseudo-holomorphic discs passing through a point of $L$ is
non-zero $($and the sum of these counts is $2^n+1)$. Therefore, for different
$n$ these tori cannot be related by symplectomorphisms.
\end{theorem}

\begin{remark}\ %
\begin{enumerate}
\item Taking the product of these tori with circles of the appropriate radius,
we also obtain infinitely many examples in $\R^{2d}$ for all $2d\ge 6$
(similarly distinguished by counts of Maslov index 2 pseudo-holomorphic
discs).
\smallskip
\item For $n=1$ our tori are most likely symplectomorphic to standard product
tori. For $n=0$ they can be shown to be symplectomorphic to the product of
a circle in $\R^2$ with the monotone Chekanov torus in $\R^4$.
\smallskip
\item Vianna's recent result concerning the existence of infinitely many 
monotone Lagrangian tori in $\CP^2$ (\cite{Vianna2}, see also \cite{GM}) should also imply a result similar to Theorem
\ref{thm:main}, by considering the preimages of these tori under the
natural projection map from the unit sphere $S^5\subset \R^6$ to $\CP^2$.
However, the construction we give here is substantially simpler.
\smallskip
\item Monotonicity plays a key role in the construction. Indeed, after arbitrarily small 
Lagrangian isotopies (not preserving monotonicity), our tori become Hamiltonian
isotopic to standard product tori.
\smallskip
\item The least elementary part of our argument is the discussion of
orientations of moduli spaces. The
reader unwilling to delve into these should be content to work with
mod 2 counts of holomorphic discs; 
the number of Maslov index~2 classes for which the algebraic count 
of discs is non-zero mod 2, 
and the number of integer points in their convex hull inside 
$\pi_2(\R^6,L)\simeq \Z^3$, are in fact
sufficient to distinguish the monotone tori we construct
for different $n$.

\end{enumerate}
\end{remark}

\subsection*{Acknowledgements}
While the methods of this paper are elementary, some of the key conceptual
ideas come from the joint work of the author with Mohammed Abouzaid and 
Ludmil Katzarkov \cite{AAK}, and from Renato Vianna's thesis \cite{Vianna}
(see \S\ref{s:motivation}).
I thank all three of them for helping shape my thoughts
on this subject. I also thank Felix Schlenk and the anonymous referees
for their careful comments.

This work was partially supported by NSF grants DMS-1406274 and
DMS-1264662, and by a fellowship from the Simons Foundation.

\section{K\"ahler reduction and Moser flow on the reduced space}
Our main object of study is the manifold 
\begin{equation}
X=\{(x,y,z,w)\in\C^4\,|\,xy=h(z,w)\},
\end{equation}
where for $n\ge 0$,
\begin{equation}\label{eq:hzw}
h(z,w)=cz^n+c^{-1}w-1,
\end{equation}
for $c\gg 1$ a constant (e.g., $c=10$).
As a complex manifold $X$ is isomorphic to $\C^3$ via projection to the
coordinates $(x,y,z)$, as $w=c(xy+1)-c^2z^n$.
We equip $X$ with the K\"ahler form
\begin{equation}\label{eq:omega_X}
\omega_X=\frac{i}{2} dz\wedge d\bar{z}+\frac{i}{2} dw\wedge d\bar{w}+
\kappa\left(\frac{i}{2} dx\wedge d\bar{x}+\frac{i}{2} dy\wedge d\bar{y}\right),
\end{equation}
where $\kappa>0$ is a small positive constant to be determined below.
We note that up to a rescaling of the $x$ and $y$ coordinates $\omega_X$ is
simply the restriction to $X$ of the standard K\"ahler form of $\C^4$.

The action of $S^1$ on $X$ by
\begin{equation}
e^{i\theta}\cdot(x,y,z,w)=(e^{i\theta}x,e^{-i\theta}y,z,w)
\end{equation}
is Hamiltonian, with moment map
\begin{equation}
\mu_X=\frac{\kappa}{2}(|x|^2-|y|^2).
\end{equation}
We will consider the reduced space
\begin{equation}
X_{red}=\mu_X^{-1}(0)/S^1.
\end{equation}
As a complex manifold, $X_{red}$ can be naturally identified with $\C^2$ via
projection to the coordinates $(z,w)$. 
Indeed, for fixed $(z,w)$ the part of the conic $xy=h(z,w)$
where $|x|=|y|$ consists of a single $S^1$-orbit; the reduced space
is therefore naturally a smooth complex manifold, even though $\mu_X^{-1}(0)$
is singular at the fixed points of the $S^1$-action, i.e.\ where 
$h(z,w)=0$ and $x=y=0$.

\begin{lemma}
The reduced K\"ahler form on $X_{red}\simeq \C^2$ is given by
\begin{equation}\label{eq:omega_red}
\omega_{red}=\frac{i}{2} dz\wedge d\bar{z}+\frac{i}{2} dw\wedge d\bar{w}+
\frac{i\kappa}{4} \frac{dh\wedge d\bar{h}}{|h|}
=\omega_0+\frac{\kappa}{2}\,dd^c(|h|).
\end{equation}
$($As expected this form is singular along the complex curve $h(z,w)=0.)$
\end{lemma} 

\proof
Given any point of $X_{red}$ where $h(z,w)\neq 0$, we choose a local
square root of $h$, and observe that a local section of the quotient map from
$\mu_X^{-1}(0)$ to $X_{red}$ is given by
setting $x=y=h(z,w)^{1/2}$. By definition, the reduced K\"ahler form 
$\omega_{red}$ agrees with the pullback of
$\omega_X$ under this local section map. Setting
$x=y=h^{1/2}$, we find that
$$
dx\wedge d\bar{x}+dy\wedge d\bar{y}=2d(h^{1/2})\wedge d(\bar{h}{}^{1/2})=
\frac{1}{2|h|}dh\wedge d\bar{h}.
$$
The first part of \eqref{eq:omega_red} follows immediately by substitution into 
\eqref{eq:omega_X}.
The second equality follows from the observation that
$$dd^c(|h|)=2i\partial\bar{\partial}(h^{1/2}\cdot \bar{h}^{1/2})=\frac{i}{2|h|} dh\wedge
d\bar{h}.$$\vskip-1.5ex
\endproof

Next we recall the following explicit form of Moser's lemma in the K\"ahler case.
\begin{lemma}\label{l:moser}
Let $\omega_0$ and $\omega_1=\omega_0+dd^c\varphi$ be two K\"ahler forms
on a complex manifold. Denote by $g_t=(1-t)g_0+tg_1$ 
the K\"ahler metric corresponding to the K\"ahler form $\omega_t=
\omega_0+t\,dd^c\varphi$ for $t\in [0,1]$, by $\xi_t=-\nabla_{g_t}(\varphi)$
the gradient of $\varphi$ with respect to $g_t$, and by 
$\psi_t$ the isotopy generated by $\xi_t$ wherever it is well-defined. Then
$\psi_t^*(\omega_t)=\omega_0$. 
Moreover, when $\omega_0=d\theta_0$ is exact, setting
$\theta_t=\theta_0+t d^c\varphi$, 
the pullback $\psi_t^*(\theta_t)$ differs from $\theta_0$ by an exact form.
\end{lemma}

\proof
The result follows from Moser's trick and the observation that
$$\omega_t(\xi_t,\cdot)=-g_t(\xi_t,J\,\cdot)=d\varphi(J\,\cdot).$$
Thus, $\iota_{\xi_t}\omega_t=-d^c\varphi$, and
$$\tfrac{d}{dt}(\psi_t^*\omega_t)=\psi_t^*(\tfrac{d}{dt}\omega_t+
L_{\xi_t}\omega_t)=\psi_t^*(dd^c\varphi+d\iota_{\xi_t}\omega_t)=0.$$
Similarly, in the exact case,
$$\tfrac{d}{dt}(\psi_t^*\theta_t)=\psi_t^*(\tfrac{d}{dt}\theta_t+
L_{\xi_t}\theta_t)=\psi_t^*(d^c\varphi+\iota_{\xi_t}(d\theta_t)+
d(\iota_{\xi_t}\theta_t))=\psi_t^*(d\iota_{\xi_t}\theta_t)$$
is exact as claimed.
\endproof

Applying this to the case at hand, we obtain:

\begin{lemma}\label{l:isotopyred}
Let $U$ be the complement of an arbitrarily small neighborhood of
$h^{-1}(0)$ inside an arbitrarily large ball in $\C^2$. Then there exists
a constant $\kappa_0>0$ $($depending on $U)$ and an isotopy 
$(\psi_\kappa)_{\kappa\in [0,\kappa_0]}$ defined on $U$,
$\psi_0=\mathrm{id}$, such that for
all $\kappa\in (0,\kappa_0)$,
$\psi_\kappa$ gives an exact symplectomorphism between $U\subset
(\C^2,\omega_0)$ and $\psi_\kappa(U)\subset (X_{red},\omega_{red})$.
\end{lemma}

\proof
Let $\Omega$ be a compact subset of $\C^2\setminus h^{-1}(0)$ whose interior
contains the closure of $U$. On $\Omega$, the function $|h|$ is smooth and 
has bounded derivatives, and the K\"ahler metric $g_\kappa$
associated to $\omega_{red}=\omega_0+\frac{\kappa}{2}dd^c(|h|)$ is bounded between fixed multiples of 
the standard metric $g_0$ for all $\kappa\in [0,1]$.
Thus, the vector field $\xi_\kappa=-\frac12\nabla_{g_\kappa}|h|$ is smooth
and has bounded norm on $\Omega$. Applying Lemma \ref{l:moser}, the isotopy
$\psi_\kappa$ generated by $\xi_\kappa$ is well-defined on $U$ for small
enough $\kappa$ and gives the desired symplectomorphisms.
\endproof

\section{Monotone tori in $X_{red}$ and $X$}\label{s:mainstuff}

\subsection{An enumerative invariant of monotone Lagrangians}\label{ss:count}
Before proceeding with our construction, we recall some basic facts
about holomorphic discs and the invariant we use to distinguish our tori.
(See also \cite{EP,Au,Vianna}.)

Let $L$ be a closed oriented spin Lagrangian
submanifold in a symplectic manifold $(M^{2d},\omega)$ equipped with a
compatible almost-complex structure $J$. When $M$ is non-compact we always
assume that $\omega$ is convex at infinity (in our case, this follows from
the properness and strict plurisubharmonicity of the K\"ahler potential).

Given a $J$-holomorphic map
 $u:(D^2,\partial D^2)\to (M,L)$, the Maslov index
$\mu([u])\in2\Z$ is the homotopy class of the loop of Lagrangian spaces given by
$TL$ along the boundary of $u$ (relative to a
trivialization of $u^*TM$). The deformation of $u$ as a $J$-holomorphic map
is governed by a Cauchy-Riemann type operator (in the integrable case, an
honest $\dbar$ operator) on the space of sections of $u^*TM$ taking
values in $u^*TL$ along the boundary. The index of this operator is
$\mathrm{ind}(\dbar)=d+\mu([u])$, and when it is surjective (i.e., $u$ is {\em
regular}) the space of pseudo-holomorphic maps is locally a smooth manifold of
this dimension. 

Assume now that $L$ is monotone, and fix a homotopy class 
$\beta\in \pi_2(M,L)$ with $\mu(\beta)=2$. We consider the moduli space
of $J$-holomorphic discs with one boundary marked point $1\in \partial D^2$,
i.e.\ the quotient
\begin{equation}
\mathcal{M}_1(L,\beta,J)=\bigl\{u:(D^2,\partial D^2)\to (M,L)\,|\,\dbar_J u=0,\ 
u_*[D^2]=\beta\bigr\}/\mathrm{Aut}(D^2,1).
\end{equation}
Since $\mu(\beta)=2$ takes the smallest possible positive value, and 
the monotonicity of $L$ guarantees that the symplectic area of discs is positively
proportional to their Maslov index, 
discs in the class $\beta$ have the smallest possible
symplectic area. Therefore, bubbling can be excluded {\em a
priori}. Moreover, all $J$-holomorphic discs in the class $\beta$ are
somewhere injective, and so a generic choice of $J$ ensures their
regularity. $\mathcal{M}_1(L,\beta,J)$ is then a smooth compact manifold
of dimension $d+\mu(\beta)-2=d$. 

Fix an orientation and a spin structure on $L$. The spin structure 
determines an orientation of
$\mathcal{M}_1(L,\beta,J)$ (cf.\ \cite{FO3book,Cho}), and
the degree of the evaluation map 
\begin{align*}
ev:\mathcal{M}_1(L,\beta,J)&\to L,\\ [u]&\mapsto u(1)
\end{align*}
is then a well-defined integer -- essentially, a signed count of
$J$-holomorphic discs in the class $\beta$ whose boundary
passes through a given point of $L$. Moreover, a generic path between two
regular almost-complex structures $J_0$ and $J_1$ determines an oriented
cobordism between $\mathcal{M}_1(L,\beta,J_0)$ and
$\mathcal{M}_1(L,\beta,J_1)$, which shows that the degree of the evaluation
map is independent of the chosen regular $J$. We denote its value by
$n(L,\beta)\in\Z$.

\begin{definition}
We call $n(L,\beta)\in\Z$ the {\em algebraic count} of pseudo-holomorphic
discs in the class $\beta$ passing through a point of $L$.
\end{definition}

By the same cobordism argument, the algebraic counts $n(L,\beta)$ are
invariant under isotopies of $L$ among monotone Lagrangian submanifolds;
and they are also invariant under simultaneous deformations of the
symplectic form on $M$ and of the Lagrangian submanifold $L$, as long as
convexity at infinity and monotonicity are preserved.
Another invariance property concerns symplectomorphisms of $M$:
if $L'=\phi(L)$ for some symplectomorphism $\phi$, then
$\mathcal{M}_1(L,\beta,J)\simeq \mathcal{M}_1(L',\phi_*\beta,\phi_*J)$, and
so (with compatible choices of orientations and spin structures) we have $n(L,\beta)=n(L',\phi_*\beta)$.

As pointed out in the introduction, the reader unwilling to deal with spin
structures and orientations of moduli spaces should be content to work with
$n(L,\beta)\ \text{mod}\ 2$.

\subsection{A monotone torus in $X_{red}$}

Let $T_{std}=\{(z,w),\ |z|=|w|=1\}$ be the standard product torus in
$(\C^2,\omega_0)$ equipped with the standard K\"ahler form and the standard
complex structure. The following is well-known (see e.g.\ \cite{Cho}; we
sketch the proof for completeness):

\begin{lemma}\label{l:Tstd}
$T_{std}$ is a monotone Lagrangian torus in $(\C^2,\omega_0)$.
There are two families of holomorphic discs of Maslov index $2$ with boundary
on $T_{std}$, which can be parametrized by the maps 
$u_\alpha:z\mapsto (z,e^{i\alpha})$ and
$v_\alpha:z\mapsto (e^{i\alpha},z)$ for $e^{i\alpha}\in S^1$. These discs
are all regular, and for a suitable choice of spin structure on $T_{std}$
the algebraic count of discs passing through a point of $T_{std}$ is $+1$ for each of
the two families.
\end{lemma}

\proof
The maps $u_\alpha$ and $v_\alpha:(D^2,\partial D^2)\to (\C^2,T_{std})$ 
obviously define holomorphic discs. To calculate their Maslov index, we note
that the pullback bundle $u_\alpha^*(T\C^2)$ can be identified with 
the direct sum of two trivial holomorphic line bundles in such a way that,
at a point $e^{i\theta}\in \partial D^2$, the pullback of 
$TT_{std}$ splits into the direct sum of the real lines
$\ell_1=e^{i\theta}\R\subset \C$ in the first factor and 
$\ell_0=\R\subset \C$ in the second factor.

Thus, the Maslov index of $u_\alpha$ is equal to the sum of the Maslov indices of
the two families of lines $\ell_1$ and $\ell_0$ in $\C$, namely $2+0=2$. 
Furthermore, the regularity of $u_\alpha$ follows from the surjectivity 
of the $\dbar$ operator for complex-valued functions on the disc with boundary 
conditions in $\ell_1$ (resp.\ $\ell_0$) (as follows e.g.\ from the
reflection principle). 
Similarly for $v_\alpha$.

To see that these are the only Maslov index 2 discs, we observe that
$\beta_1=[u_\alpha]$ and $\beta_2=[v_\alpha]$ generate $\pi_2(\C^2,T_{std})\simeq \pi_1(T_{std})=\Z^2$,
so by linearity the Maslov index of a disc with boundary on
$T_{std}$ is equal to twice its algebraic intersection number with the
union of the coordinate axes. For holomorphic discs, positivity of
intersection implies that a Maslov index 2 disc in $(\C^2,T_{std})$
intersects only one of the two coordinate axes $z=0$ and $w=0$,
transversely, and at a single point. 

If for example the holomorphic disc
$u:(D^2,\partial D^2)\to (\C^2,T_{std})$ is disjoint
from the line $w=0$, then applying the maximum principle to the
projection to the $w$ coordinate, we find that $w\circ u:(D^2,\partial D^2)\to
(\C^*,S^1)$ must take some constant value $e^{i\alpha}$. Meanwhile, the
projection to the $z$ coordinate has a single zero of order 1, which means that
$z\circ u:(D^2,\partial D^2)\to (\C, S^1)$ is a biholomorphism from the unit
disc to itself, i.e.\ the identity map up to reparametrization. Thus $u$ is
equivalent to $u_\alpha$ up to reparametrization. Similarly for the other
case where the disc is disjoint from $z=0$ and intersects $w=0$ once.

Finally, the moduli space $\mathcal{M}_1(L,\beta_1,J_0)$ consists of 
reparametrizations of the discs $u_\alpha$, e.g.\ the maps $z\mapsto
(e^{i\beta}z,e^{i\alpha})$ for $(e^{i\beta},e^{i\alpha})\in S^1\times S^1$.
Thus $\mathcal{M}_1(L,\beta_1,J_0)\simeq T^2$, and the evaluation map to
$T_{std}$ is a diffeomorphism; choosing the ``standard'' spin structure
ensures that this diffeomorphism is orientation-preserving \cite{Cho},
hence $n(L,\beta_1)=+1$. Similarly for the other class $\beta_2$.
\endproof

Next we observe that $T_{std}$ lies away from the complex curve
\begin{equation}
C=h^{-1}(0)=\{(z,w)\in\C^2\,|\,cz^n+c^{-1}w-1=0\},
\end{equation}
and that the disc $u_\alpha$ intersects $C$ transversely at $n$ distinct
points, where the $z$ coordinate takes the values $$z_k=e^{2\pi i k/n} c^{-1/n} (1-c^{-1}e^{i\alpha})^{1/n},$$ while
$v_\alpha$ is disjoint from $C$.

The regularity of the discs $u_\alpha$ and $v_\alpha$ implies that they
deform smoothly under small isotopies of $T_{std}$. Thus, for small enough
values of the constant $\kappa$, denoting by $\psi_\kappa$ the isotopy
constructed in Lemma \ref{l:isotopyred}, the Lagrangian torus 
$$T_{red}=\psi_\kappa(T_{std})$$ in $(X_{red},\omega_{red})$ again bounds two
families of Maslov index 2 holomorphic discs $u'_\alpha$ and
$v'_\alpha$, representing the homotopy classes
$\beta'_1=(\psi_\kappa)_*(\beta_1)$ and $\beta'_2=(\psi_\kappa)_*(\beta_2)$.
We obtain:

\begin{lemma}\label{l:Tred}
For $\kappa>0$ small enough, $(X_{red},\omega_{red})$ contains a monotone Lagrangian torus
$T_{red}$, disjoint from $C=h^{-1}(0)$,
which bounds exactly two families
of Maslov index $2$ holomorphic discs, representing classes $\beta'_1,\beta'_2$
that span $\pi_2(X_{red},T_{red})\simeq \Z^2$. 
These discs are all regular, and for a
suitable spin structure their algebraic counts are
$n(T_{red},\beta'_1)=n(T_{red},\beta'_2)=+1$. Moreover, the discs in the
class $\beta'_1$ intersect $C$ transversely in $n$ distinct
points, while those in the class $\beta'_2$ are disjoint from $C$.
\end{lemma}

\begin{remark} While $\omega_{red}$ is singular along $C$, it can still be
integrated over a disc that intersects $C$ transversely, so the notion of
monotonicity still makes sense. In fact, symplectic area can also be defined 
as the integral of the Liouville form
$$\theta_{red}=d^c(\tfrac14 |z|^2+\tfrac14 |w|^2+\tfrac{\kappa}{2}|h|)$$ along
the boundary of a disc.
Perhaps even better, we can modify $\omega_{red}$ in a neighborhood of $C$ 
(disjoint from $T_{red}$) by a small exact deformation so as to cure its
lack of smoothness; this can be achieved simply by replacing $|h|$ by a 
smooth function $\rho(|h|)$ in the expression for the K\"ahler potential
(taking $\rho:[0,\infty)\to [0,\infty)$ to be any smooth, convex function
which agrees with identity outside of $[0,\epsilon]$ and has vanishing odd derivatives at the
origin). This modification does not affect the properties of the isotopy
$\psi_\kappa$ away from $C$, nor the symplectic areas of holomorphic
discs.
\end{remark}

\proof[Proof of Lemma \ref{l:Tred}]
The existence and regularity for small $\kappa$ of the two families of
holomorphic discs $u'_\alpha$ and $v'_\alpha$ with boundary
on $T_{red}=\psi_\kappa(T_{std})$ representing the classes $\beta'_1$ and 
$\beta'_2$, obtained as smooth deformations of the discs $u_\alpha$ and
$v_\alpha$ under the isotopy, is a direct consequence of the regularity of 
the latter discs. 

Since the isotopy is exact ($\psi_\kappa^*(\theta_{red})$
agrees with the standard Liouville form $\theta_0$ up to an exact term),
the symplectic areas of the discs are preserved, which proves the
monotonicity of $T_{red}$. Moreover, Gromov compactness implies that $T_{red}$ does 
not bound any other Maslov index 2 holomorphic discs: if such discs
existed for arbitrarily small $\kappa$, taking the limit of a subsequence
with $\kappa\to 0$ would yield a contradiction. 

Finally, because the discs
$u_\alpha$ and $v_\alpha$ deform smoothly under the isotopy of $T_{std}$ to
$T_{red}$, for small $\kappa$ the discs $u'_\alpha$ and $v'_\alpha$ continue to intersect $C$
transversely, and the algebraic counts remain unchanged (in fact the
evaluation maps $ev:\mathcal{M}_1(T_{red},\beta'_i,J_0)\to T_{red}$ remain diffeomorphisms).
\endproof

\subsection{A monotone torus in $X$}
From now on we fix the value of the constant $\kappa>0$ so that the
conclusion of Lemma \ref{l:Tred} holds. We then construct a Lagrangian torus
$T$ in $(X,\omega_X)$ by lifting $T_{red}$ to $\mu_X^{-1}(0)$:

\begin{definition}
We denote by $T$ the preimage of $T_{red}$ under the projection map from
$\mu_X^{-1}(0)\subset X$ to $X_{red}$, i.e.\
\begin{equation}
T=\{(x,y,z,w)\in X\,|\,(z,w)\in T_{red}\text{ and }|x|=|y|\}.
\end{equation}
\end{definition}

\noindent We also denote by $\pi:X\to X_{red}$ the projection to the $(z,w)$
coordinates, 
\begin{equation}
\pi(x,y,z,w)=(z,w).
\end{equation}

\begin{lemma}\label{l:Tmonotone}
$T$ is a monotone Lagrangian torus in $(X,\omega_X)$.
\end{lemma}

Conceptually, this follows from the observation that $T$ is the image of $T_{red}$
under the monotone Lagrangian correspondence between $X_{red}$ and $X$
induced by $\mu_{X}^{-1}(0)$. A more elementary argument is as follows.

\proof
Since the restriction of $\omega_X$ to $\mu_X^{-1}(0)$ agrees with the
pullback of $\omega_{red}$ via the projection map $\pi$, $\omega_{X}|_T$ is the
pullback of $\omega_{red}|_{T_{red}}$ under the projection from $T\subset
\mu_X^{-1}(0)$ to
$T_{red}\subset X_{red}$, i.e.\ it vanishes, and $T$ is Lagrangian.

Let $u:(D^2,\partial D^2)\to (X,T)$ be a disc with boundary on $T$ (not
necessarily holomorphic), and denote by $\gamma:S^1\to T$ its boundary loop. 
Perturbing $u$ if necessary, we can assume that it
avoids the fixed point set $F=\{x=y=0\}$ (which has real codimension~4).
In terms of the Liouville form
\begin{equation}
\theta_X=d^c(\tfrac14 |z|^2+\tfrac14 |w|^2+\tfrac{\kappa}{4} |x|^2
+\tfrac{\kappa}{4}|y|^2),
\end{equation}
the symplectic area of $u$ is given by the integral of $\theta_X$ along the
boundary loop~$\gamma$.
However, along $\mu_X^{-1}(0)$ we have $|x|^2=|y|^2=|h|$, and 
$|x|^2+|y|^2$ achieves its fiberwise minimum so its derivative
vanishes in all directions tangent to the fibers of $\pi$.
Therefore, at every point of $\mu_X^{-1}(0)$ the 1-form
$\theta_X$ coincides with
$$\pi^*\theta_{red}=d^c(\tfrac14|z|^2+\tfrac14|w|^2+\tfrac{\kappa}{2}|h|).$$
Denoting by $u_{red}=\pi\circ u:(D^2,\partial D^2)\to (X_{red},T_{red})$ and
$\gamma_{red}=\pi\circ\gamma:S^1\to T_{red}$ the projections of $u$ and
$\gamma$, we conclude that
\begin{equation}\label{eq:equalareas}
\int_{D^2} u^*\omega_X=\int_{S^1} \gamma^*\theta_X=
\int_{S^1} \gamma^*(\pi^*\theta_{red})=\int_{S^1}
\gamma_{red}^*(\theta_{red})=\int_{D^2} u_{red}^*(\omega_{red}),
\end{equation}
i.e.\ the disc $u$ and its projection $u_{red}$ have the same symplectic
areas. Meanwhile, away from the fixed point locus $F$, denote by
\begin{equation}
\mathcal{L}_\R=\R\cdot (ix,-iy,0,0)\quad \text{and}\quad 
\mathcal{L}=\C\cdot (ix,-iy,0,0)
\end{equation}
the real and complex
spans of the vector field generating the $S^1$-action. Then
$\mathcal{L}$ is a trivial holomorphic subbundle of $TX$, and
$TX/\mathcal{L}\simeq \pi^*TX_{red}$, i.e.\ away from $F$ we have a
short exact sequence of holomorphic vector bundles
\begin{equation}\label{eq:seshol}
0\longrightarrow \mathcal{L}\longrightarrow 
TX\stackrel{d\pi}{\longrightarrow} \pi^* TX_{red}\longrightarrow 0.
\end{equation}
Along $T$, we have a similar short exact sequence of real subbundles,
\begin{equation}\label{eq:sesreal}
0\longrightarrow \mathcal{L}_\R\longrightarrow TT
\stackrel{d\pi}{\longrightarrow} \pi^*TT_{red}\longrightarrow 0.
\end{equation}
Since the trivial subbundles $(u^*\mathcal{L},\gamma^*\mathcal{L}_\R)$ do
not contribute to the Maslov index, $\mu([u])$ can be computed
by considering the quotient bundles $(u^*(TX/\mathcal{L}),\gamma^*(TT/\mathcal{L}_\R))
\simeq (u_{red}^*(TX_{red}),\gamma_{red}^*(TT_{red}))$. In other terms,
\begin{equation}\label{eq:equalmaslov}
\mu([u])=\mu([u_{red}]).
\end{equation}
Comparing \eqref{eq:equalareas} and \eqref{eq:equalmaslov},
we find that the proportionality between Maslov index and symplectic area for discs
in $X_{red}$ with boundary on $T_{red}$ implies the same proportionality for
discs in $X$ with boundary on $T$.
\endproof

\begin{lemma}\label{l:liftdiscs}
The projection 
$u_{red}=\pi\circ u:(D^2,\partial D^2)\to (X_{red},T_{red})$ of a
holomorphic disc $u:(D^2,\partial D^2)\to (X,T)$ is a
holomorphic disc, and $\mu([u_{red}])=\mu([u])$.

Conversely, let $u_{red}:(D^2,\partial D^2)\to
(X_{red},T_{red})$ be a holomorphic disc that 
intersects $C=h^{-1}(0)$ transversely in $k$ 
points, and fix a point $p_0\in T$ such that $\pi(p_0)=u_{red}(1)$.
Then there are exactly $2^k$
holomorphic discs $u:(D^2,\partial D^2)\to (X,T)$ such that
$\pi\circ u=u_{red}$ and $u(1)=p_0$. Moreover, if $u_{red}$ is regular
then all these discs are regular.
\end{lemma}

\proof
The first statement follows immediately from the holomorphicity of $\pi$ and
the Maslov index calculation in the proof of Lemma \ref{l:Tmonotone}
(equation \eqref{eq:equalmaslov}).

For the second part, let $u_{red}$ be a holomorphic disc in $X_{red}$ that
intersects $C$ transversely, with $u_{red}^{-1}(C)=\{t_1,\dots,t_k\}\subset
D^2$, and let $u$ be a lift of $u_{red}$ to a disc in $X$ with boundary on
$T$.  Along the holomorphic disc $u$, the product $xy=h(z,w)$ has
simple zeroes at $t_1,\dots,t_k$, i.e.\ $u$ intersects
$\pi^{-1}(C)=\{x=0\}\cup \{y=0\}$ transversely at the $k$ points
$u(t_1),\dots,u(t_k)$.
The quotient $q=x/y$ then defines a meromorphic function on the disc, which
has either a simple zero or a simple pole at each of $t_1,\dots,t_k$, and
no other zeroes or poles. Moreover, on the boundary we have $|x|=|y|$, so
$q$ maps the unit circle to itself. 

Given any function $\varepsilon:\{1,\dots,k\}\to \{\pm 1\}$, set
\begin{equation}
\vartheta_\varepsilon(z)=\prod_{j=1}^k
\left(\frac{z-t_j}{1-\overline{t_j}z}\right)^{\varepsilon(j)},
\end{equation}
which is a meromorphic function on the unit disc, mapping the unit circle
to itself, and with simple zeroes (resp.\ poles) at all $t_j$ such that $\varepsilon(j)=+1$
(resp.\ $-1$).

Thus, choosing $\varepsilon(j)=\mathrm{ord}_{t_j}(q)$ according to the poles
and zeroes of $q=x/y$ along the disc $u$, we find that $\vartheta_\varepsilon$ and
$q$ have the same zeroes and poles on the unit disc, and their ratio defines
a nowhere vanishing holomorphic function on the unit disc, taking values in
the unit circle at the boundary. By the maximum principle this function is
constant, i.e.\ there exists $e^{i\theta}\in S^1$ such that
$q=e^{i\theta}\vartheta_\varepsilon$. 

By construction the holomorphic functions $(h\circ u_{red})\vartheta_\varepsilon^{\pm 1}$ 
only have double zeroes, and so we can choose square roots 
$$\zeta_{\pm}= \left((h\circ u_{red})\,\vartheta_{\varepsilon}^{\pm 1}\right)^{1/2},$$
with $\zeta_+/\zeta_-=\vartheta_{\varepsilon}$ and
$\zeta_+\zeta_-=h\circ u_{red}$. We obtain that
along the disc $u$ the coordinates $x$ and $y$ are given by
$$x=e^{i\theta/2} \zeta_+\quad \text{and}\quad y=e^{-i\theta/2} \zeta_-,$$
for some $e^{i\theta/2}\in S^1$.
Conversely, these formulas determine holomorphic lifts of $u_{red}$ for all
$\varepsilon:\{1,\dots,k\}\to \{\pm 1\}$ and for all $e^{i\theta/2}\in S^1$,
and the condition that $u(1)=p_0$ determines the normalization factor
$e^{i\theta/2}$ uniquely for given $\varepsilon$. Hence there are $2^k$
lifts of $u_{red}$ as claimed, determined by the choice of whether $x$ or
$y$ vanishes at each point where $u_{red}$ intersects $C$.

Finally, we note that none of the lifts $u$ pass through the fixed point
locus of the $S^1$-action (since $x$ and $y$ do not vanish simultaneously).
Thus, pulling back the exact sequences \eqref{eq:seshol}
and \eqref{eq:sesreal} along $u$, we find that the holomorphic vector bundle
$u^*TX$ admits a trivial holomorphic line subbundle $u^*\mathcal{L}$, with
a trivial real subbundle at the boundary $u_{|S^1}^*\mathcal{L}_\R$. Since
the $\dbar$ operator for complex-valued functions on the disc with the
trivial real  boundary condition $\R\subset\C$ on the unit circle is 
surjective, the surjectivity of the $\dbar$ operator on sections of
$u^*TX$ with boundary conditions $u_{|S^1}^*(TT)$ is equivalent to that of
the $\dbar$ operator on the quotient bundle 
$u^*TX/u^*\mathcal{L}\simeq u_{red}^*TX_{red}$ with boundary conditions
$u_{|S^1}^*(TT)/u_{|S^1}^*(\mathcal{L}_\R)\simeq u_{red|S^1}^*(TT_{red})$.
Thus, the regularity of $u$ is equivalent to that of $u_{red}$ as claimed.
\endproof

\begin{corollary}\label{cor:TX}
There are $n+2$ distinct Maslov index\/ $2$ classes in $\pi_2(X,T)$ for which
the algebraic count of pseudo-holomorphic discs is non-zero, and for a
suitable choice of spin structure the sum of these counts is $2^n+1$.
\end{corollary}

\proof
By Lemma \ref{l:liftdiscs}, the holomorphic discs of Maslov index 2 bounded 
by $T$ are lifts of those bounded by $T_{red}$ in $X_{red}$, which are
determined by Lemma \ref{l:Tred}. 

The discs representing the class
$\beta'_2\in \pi_2(X_{red},T_{red})$ are disjoint from $C$, hence they admit
a unique lift up to the $S^1$-action. Denoting by
$\hat{\beta}_2\in \pi_2(X,T)$ the class of these lifts, the moduli space 
$\mathcal{M}_1(T,\hat{\beta}_2,J_0)$ is an $S^1$-bundle over
$\mathcal{M}_1(T_{red},\beta'_2,J_0)$, and the evaluation map to $T$ is
equivariant with respect to the $S^1$-action; thus the evaluation map
$ev:\mathcal{M}_1(T,\hat{\beta}_2,J_0)\to T$ is again a diffeomorphism,
and its degree is $\pm 1$.

Meanwhile, the discs representing the class $\beta'_1\in \pi_2(X_{red},T_{red})$
intersect $C$ transversely in $n$ points (cf.\ Lemma \ref{l:Tred}), so by
Lemma \ref{l:liftdiscs} they can be lifted in $2^n$ different ways up to the
$S^1$-action. Observe that elements of $\pi_2(X,T)\simeq \Z^3$ are determined by
their intersection numbers with the three hypersurfaces $x=0$, $z=0$, and
$w=0$. Thus, the lifts live in $n+1$ different classes
$\hat{\beta}_{1,\ell}\in \pi_2(X,T)$, $\ell=0,\dots,n$, depending on the
intersection number of the lifted disc with the hypersurface $x=0$;
each value of $\ell$ is achieved by $\binom{n}{\ell}$ of the $2^n$ lifts.
The moduli space $\mathcal{M}_1(T,\hat{\beta}_{1,\ell},J_0)$ then
projects to $\mathcal{M}_1(T_{red},\beta'_1,J_0)$ with fiber a union of
$\binom{n}{\ell}$ circles. The evaluation map
$ev:\mathcal{M}_1(T,\hat{\beta}_{1,\ell},J_0)\to T$ is thus an unramified
$\binom{n}{\ell}$-sheeted covering. 

To determine the orientations, we briefly recall the construction in
\cite[Chapter~8]{FO3book} (see also \cite[Prop.~5.2]{Cho} for a simpler
presentation that suffices for the case at hand). 
A spin structure on $T$
determines a trivialization of its tangent bundle along the boundary of a
holomorphic disc $u$. Using this trivialization, the $\dbar$ 
operator can be deformed to the direct sum of a complex linear operator and a $\dbar$
operator for sections of a trivialized complex vector bundle with trivial
real boundary condition (namely, the tangent bundles to $X$ and $T$ along
the boundary of $u$, with the trivialization determined by the spin
structure). Since the kernel of the latter operator can be identified with
the tangent space to $T$ at the marked point, an orientation of $T$ then 
determines an orientation of the tangent space to the moduli space at $u$.

In our case, we choose the spin structure
on $T$ to be standard along the orbits of the $S^1$-action and consistent
under the splitting \eqref{eq:sesreal} with that previously chosen on $T_{red}$. Thus, the preferred trivialization of $TT$
along the boundary of a holomorphic disc $u$ agrees with that induced via
\eqref{eq:sesreal} by
the trivialization of $TT_{red}$ along the boundary of $u_{red}=\pi\circ u$
and the natural trivialization of the trivial line bundle $\mathcal{L}_\R$.
The orientation at $u$ of the moduli space of holomorphic discs in $(X,T)$ 
then agrees with that induced by the orientation at $u_{red}$ of the moduli
space of holomorphic discs in $(X_{red},T_{red})$ and the chosen orientation
of the orbits of the $S^1$-action. With this understood, the
orientation-preserving nature of the
evaluation maps for discs in $(X_{red},T_{red})$ implies that the evaluation
maps for discs in $(X,T)$ are also orientation-preserving, i.e.\ the degrees
are positive.
\endproof

\noindent
(For the reader working mod 2, we note that the odd values of $n(T,\beta)$
are achieved for $\hat{\beta}_2$ and those $\hat{\beta}_{1,\ell}$ for which
$\binom{n}{\ell}$ is odd, including the extremal cases $\hat{\beta}_{1,0}$ and
$\hat{\beta}_{1,n}$.)

\section{Proof of Theorem \ref{thm:main}}

In light of Corollary \ref{cor:TX} and the invariance properties of the algebraic
counts $n(T,\beta)$, the only thing that remains to be
done is to construct an isotopy between the K\"ahler form $\omega_X$ 
on (a bounded subset of) $X\simeq \C^3$ and the standard K\"ahler form. We will again rely on
Moser's trick (Lemma \ref{l:moser}).
We denote by 
\begin{equation}
\Phi_1=\tfrac{\kappa}{4}|x|^2+\tfrac{\kappa}{4}|y|^2+\tfrac14
|z|^2
\end{equation} the K\"ahler potential for the standard (up to rescaling) K\"ahler
form on $\C^3$,
$$\omega_1=dd^c\Phi_1=\tfrac{i}{2} dz\wedge d\bar{z}+
\kappa(\tfrac{i}{2}dx\wedge d\bar{x}+\tfrac{i}{2}dy\wedge d\bar{y}).$$
The K\"ahler potential for $\omega_X$ is
$$\Phi_X=\Phi_1+\tfrac14 |w|^2,$$
where we recall that $w$ is determined as a function of the coordinates
$(x,y,z)$ by 
\begin{equation}
w=c(xy+1)-c^2z^n.
\end{equation}

The estimate that ensures the existence of the Moser flow is the following:
\begin{lemma}\label{l:barrier}
Given any bounded subset $B\subset \C^3$,
there exist positive constants $C$ and $M$ such that the real-valued
function $\varphi=C\Phi_1-\Phi_X$ is bounded above by $M$ on $B$, and
the connected component $\Omega$ of $\varphi^{-1}((-\infty,M])$ which contains $B$
is compact.
\end{lemma}

\proof We equip $\C^3$ with the Euclidean metric for which the positive
definite quadratic form $\Phi_1$ is the square of the distance to the origin
(i.e., a rescaling of the usual metric). 

Let $R>0$ be such that $B$ is contained within the ball
$B(0,R)$ of radius $R$ (for this metric), denote by $K$ the supremum of
$\frac14|w|^2/\Phi_1$ in $B(0,2R)\setminus B(0,R)$, and set $C=2K+1$. 
Then in $B(0,2R)\setminus B(0,R)$ we have
$$K\Phi_1\le \varphi=(C-1)\Phi_1-\tfrac14 |w|^2 \le 2K\Phi_1,$$
and the upper bound continues to hold inside $B(0,R)$.

Then inside $B(0,R)$
we have $\varphi\le 2K\Phi_1\le 2KR^2$, while
in $B(0,2R)\setminus B(0,\sqrt{3}R)$ we have $3KR^2\le K\Phi_1\le \varphi$.
Thus, setting $M=\frac52 KR^2$, there is a connected component
$\Omega$ of $\varphi^{-1}((-\infty,M])$ for which $B(0,R)\subset \Omega
\subset B(0,\sqrt{3}R)$. 
\endproof

Choosing $B$ to be a polydisc in $\C^3$ large enough to contain $T$, and
taking $C$ as in Lemma \ref{l:barrier}, 
we now apply Lemma \ref{l:moser} to the K\"ahler forms $\omega_X$ and
$C\omega_1$, to construct an exact isotopy $\psi_t$ such that
$\psi_1^*(C\omega_1)=\omega_X$. 
Because the isotopy is generated by the negative gradient of
$\varphi=C\Phi_1-\Phi_X$ (with respect to a varying family of K\"ahler
metrics), the values of $\varphi$ decrease along the flow. Thus, the
compact subset $\Omega\supset B$ constructed in Lemma \ref{l:barrier} is
preserved, so the isotopy is well-defined everywhere in it, and in
particular in $B$. 

Since the isotopy is exact, $\psi_t(T)$ is a monotone Lagrangian torus
in $\C^3$ equipped with the K\"ahler form
$\omega_t=Ct\omega_1+(1-t)\omega_X$,
and the algebraic counts of Maslov index 2 holomorphic discs remain 
constant along the isotopy. For $t=1$ we obtain a monotone Lagrangian
torus in $(\C^3,C\omega_1)$ with the desired properties. Rescaling the
coordinate axes by suitable constant factors, we obtain a monotone
Lagrangian torus in $\C^3$ equipped with the standard K\"ahler form, and
by further rescaling we obtain tori with arbitrary
monotonicity constants and the same algebraic counts of pseudo-holomorphic
discs.

\section{Comments on the construction} \label{s:motivation}

Our construction is inspired by
ideas from mirror symmetry, and more precisely the Strominger-Yau-Zaslow (SYZ)
conjecture, whereby the mirror of a given K\"ahler manifold is
constructed geometrically from a Lagrangian torus fibration
on the complement of a complex hypersurface. The numbers of Maslov
index 2 discs bounded by the fibers exhibit discontinuities across a set of
{\em walls} which separate the fibration into {\em chambers}, each
with its own enumerative behavior; each chamber corresponds to a
distinguished coordinate chart on the mirror 
(cf.\ \cite[\S 2]{AAK} and \cite{Au}).

In a given Lagrangian fibration, the vast majority of fibers are
not monotone, and the counts of Maslov index 2 discs are not invariant under
Hamiltonian isotopies. However, by deforming the fibration suitably it is 
often possible to arrange the existence of a monotone fiber in any 
given chamber. For example, the complement of a smooth cubic in $\CP^2$
admits a Lagrangian torus fibration with 3 singular fibers and infinitely
many chambers; Vianna's constructions in \cite{Vianna,Vianna2} can be understood as
modifying the fibration to place the monotone fiber in a prescribed chamber.

The construction of Theorem \ref{thm:main} relies on the fact that $\C^3$
can be presented as a conic bundle $\{xy=h(z,w)\}$ over $\C^2$ with a
discriminant curve $h^{-1}(0)\subset \C^2$ of arbitrarily large degree.
SYZ mirror symmetry for conic bundles over toric varieties has been 
studied in detail in \cite{AAK}, where it was shown that the chamber
structure is governed by the tropical geometry of $h^{-1}(0)$
(or, in more classical terms, by the various manners in which product
tori in $\C^2$ can be linked with $h^{-1}(0)$).
Thus, by increasing the degree of $h$ we can exhibit Lagrangian torus
fibrations on open dense subsets of $\C^3$ (namely, those points where
$z$ and $w$ are nonzero) with arbitrarily many chambers. Choosing the
coefficients of $h$ suitably ensures
the existence of monotone fibers in the ``most interesting'' chamber.
(In fact, choosing $h$ to be analytic rather than algebraic one could 
obtain a single fibration with infinitely many chambers, with monotone
representatives corresponding to all the values of $n$ in our main
construction at once.)

Another perspective on the construction comes from singularity theory:
projecting the conic bundle $X\simeq \C^3$ to the coordinate $w$ presents
it as an unfolding of the $A_{n-1}$ singularity $xy=cz^n$. 
The $A_{n-1}$ Milnor fiber contains non-displaceable monotone Lagrangian tori
(cf.\ \cite[Corollary 9.1]{AAK} and \cite{LM}). The examples of Theorem
\ref{thm:main} can be obtained by transporting these tori along a circle 
in the $w$ coordinate; even though the unfolding makes the ambient manifold
contractible and the tori displaceable, the distinctive enumerative features
of the tori in the fibers persist.

\end{document}